\documentclass[12pt]{article}

\usepackage[utf8]{inputenc}

\usepackage[english]{babel}

\usepackage[shortlabels]{enumitem}
\usepackage[T1]{fontenc}

\usepackage{tikz}
\usepackage{verbatim}

\usepackage{amsmath,amsfonts,amssymb,amsthm}
\usepackage{nccmath}
\usepackage{mathrsfs,indentfirst, latexsym}
\usepackage{amscd}
\usepackage{eucal}
\usepackage{fancyhdr}

\numberwithin{equation}{section}

\newtheorem{theorem}{Theorem}
\newtheorem{lemma}[theorem]{Lemma}
\newtheorem{proposition}[theorem]{Proposition}
\newtheorem{corollary}[theorem]{Corollary}

\theoremstyle{definition}

\newtheorem{question}[theorem]{Question}

\newtheorem{remark}[theorem]{Remark}

\newcommand{\N}{\mathbb{N}} 
\newcommand{\Z}{\mathbb{Z}} 

\newcommand{\A}{\mathcal{A}}

\newcommand{\PS}{\mathcal{PS}}

\newcommand{\abs}[1]{\left\lvert#1\right\rvert}

\newcommand{\norm}[1]{\lVert#1\rVert}

\begin{document}
\title{Frequent hypercyclicity and piecewise syndetic recurrence sets}

\author{Yunied Puig}

\date{}
\maketitle

\begin{abstract}
Motivated by a question posed by Sophie Grivaux concerning the regularity of the orbits of frequently hypercylic operators, we show the following: for any operator $T$ on a separable metrizable and complete topological vector space $X$ which is both frequently hypercyclic and piecewise syndetic hypercyclic, the lower density and upper Banach density of the recurrence set $\{n\geq 1: T^n x\in U\}$   are different, for any hypercyclic vector $x\in X$ for $T$, and a certain collection of non-empty open sets $U\subseteq X$. As an immediate consequence we got a sufficient condition for a chaotic operator to be non frequently hypercyclic.

  \end{abstract}

\let\thefootnote\relax\footnote{Department of Mathematics, University of California at Riverside, 900 University Ave., Riverside, CA 92521, USA. e-mail:puigdedios@gmail.com}
\let\thefootnote\relax\footnote{This paper has been initiated while the author was a postdoc at Ben-Gurion University, Israel. Research supported by the Israel Science Foundation, grants 1/12 and 626/14.}
\section{Introduction}

The main motivation of this paper is the following question concerning the regularity of the orbits of frequently hypercyclic operators. This was posed by S. Grivaux in an earlier version of \cite{Gr} from 2013.

Does there exist a frequently hypercyclic operator $T$ on a separable Banach space $X$ whose recurrence set $N(x, U)=\{n\geq 1: T^n x\in U\}$, has positive and different lower density and upper density, for some frequently hypercyclic vector $x\in X$ and some non-empty open set $U\subseteq X$?

 During the completion of this manuscript, we have been informed by S. Grivaux about her success in solving this question in the positive by showing that this is indeed the case for operators satisfying a particular form of the Frequent Hypercyclicity Criterion. The reader can find all the details in \cite{Gr1}.

In connection with Grivaux's question, we could wonder how the lower density and upper Banach density of recurrence sets of frequently hypercyclic operators are related. So we get the following:

\begin{question}
\label{question_oya}
Does there exist a frequently hypercyclic operator $T$ on a separable Banach space $X$ whose recurrence set $N(x, U)$ has positive and different lower density and upper Banach density, for any frequently hypercyclic vector $x\in X$ and some non-empty open set $U\subseteq X$?
\end{question}

Our aim is to give a positive answer to Question \ref{question_oya}\footnote{Please note that Theorem \ref{theorem.irregular.BD} is cited in \cite{Gr1}. The citation in \cite{Gr1} corresponds to an earlier version of this manuscript where we believed that any frequently hypercyclic operator was piecewise syndetic hypercyclic. As showed in \cite{Pu3}, there exists a frequently hypercyclic operator which is not piecewise syndetic hypercyclic, so the assumption $T$ be piecewise syndetic hypercyclic should be considered in the citation of \cite{Gr1}.}.

  \begin{theorem} 
	\label{theorem.irregular.BD}
  For any operator $T$ on a separable metrizable and complete topological vector space $X$ being both frequently hypercyclic and piecewise syndetic hypercyclic, any hypercyclic vector $x$ and any non-empty open set $U\subseteq X$ satisfying that unions of finitely many iterates of the set $U$ are not dense in $X$, we have that the lower density and upper Banach density of $N(x,U)$ are different. 
  \end{theorem}
	
	Surprisingly, the proof of Theorem \ref{theorem.irregular.BD} is combinatorial number theory in essence. The point here is to look at the piecewise syndeticity level of recurrence sets of those operators which are both frequently hypercyclic and piecewise syndetic hypercyclic. Frequently hypercyclic operators have been object of systematical study in the last ten years. Yet, it is not well understood whether its  recurrence sets satisfy other properties rather than be of positive lower density. We remark that piecewise syndetic hypercyclicity and hence piecewise syndeticity level of recurrence sets of hypercyclic operators have not been investigated before.
	
	Since any chaotic operator is piecewise syndetic hypercyclic as remarked below, as a consequence of Theorem \ref{theorem.irregular.BD} we immediately got a sufficient condition for a chaotic operator to be non frequently hypercyclic. Moreover, this shows an intimate relationship between two of the most important problems in linear dynamics, i.e. the regularity of the orbits of frequently hypercyclic operators problem and the existence of chaotic operators which are not frequently hypercyclic.
	
	\begin{corollary}
	Let $T$ be a chaotic operator on a separable $F$-space $X$ such that $N(x, U)$ has equal lower density and upper Banach density for some hypercyclic vector $x$ and some non-empty open set $U$ of $X$ satisfying that unions of finitely many iterates of the set $U$ are not dense in $X$, then $T$ is not frequently hypercyclic.
	\end{corollary}
	
	\section{Preliminaries and notation}
	From now on, if $a, b\in \N$, the interval $[a,b]$ stands for $[a, b]\cap \N$ and $\abs{A}$ for the cardinality of the set $A$. 
	Let $\A$ be a \emph{Furstenberg family} of subsets of $\Z_+$ (i.e., $\emptyset\notin \A$ and for any $A\in \A$, if $A\subset B$ then $B\in \A$). From now on, $\A$ will denote a Furstenberg family of subsets of $\Z_+$. We will be interested in the following Furstenberg families:

\[\underline{\mathcal{D}}=\{A\subseteq \Z_+:\underline{d}(A)=\liminf_{n} \abs{A\cap [0, n]}/(n+1)>0\};\]

\[ \overline{\mathcal{D}}=\{A\subseteq \Z_+: \overline{d}(A)=\limsup_{n} \abs{A\cap [0, n]}/(n+1)>0\};\]

 \[\overline{\mathcal{BD}}=\{A\subseteq \Z_+: \overline{Bd}(A)=\lim_n \limsup_m \abs{A\cap [m, m+n]}/(n+1)>0\};\] commonly known as the family of sets of positive \emph{lower density}, positive \emph{upper density} and positive \emph{upper Banach density} respectively.

 According to \cite{GrToTo}, $\overline{Bd}(A)$ it is known to coincide with 
\begin{equation}
\label{condition_etna}
\sup\{\alpha\geq 0: \forall m\in \N, \exists n, r\in \N\cup \{0\}, r\geq m \quad \abs{A\cap [n, n+r)}\geq \alpha r\}.
\end{equation}

An infinite subset $A$ of $\N$ is said to be \emph{syndetic} if $\cup_{t=1}^b A-t=\N$ for some $b\geq 1$.

 An infinite subset $A$ of $\N$ is said to be \emph{thick} if it contains arbitrarily long intervals, i.e. for any $L>0$ there exists $m\geq 1$ such that the set $\{m\geq1: [m, m+L)\subset A\}$ is non-empty. It is well-known that $A$ is thick if and only if $ \overline{Bd}(A)=1$.

An important Furstenberg family for our purpose is the one of piecewise syndetic sets.
An infinite subset $A$ of $\Z_+$ is said to be \emph{b-piecewise syndetic} ($A\in b-\PS$, for short) with $b\geq 1$, if $\cup_{t=1}^b (A-t)$ is thick.

A set $A\subseteq \Z_+$ is said to be \emph{piecewise syndetic} ($A\in \PS$, for short) if it is $b$-piecewise syndetic for some $b\geq 1$. 
It is clear that a set $A$ is piecewise syndetic if and only if there exists $b\geq 1$ such that $\overline{Bd}(\cup_{t=1}^b A-t)=1$.

 Every piecewise syndetic set has clearly positive upper Banach density, however the converse is not true. Moreover, it is known that we can have sets of positive lower density which are not piecewise syndetic. For example, the complement of the set 
  $E=\cup_{j=n}^\infty\cup_{m\geq 1}[2^jm, 2^jm+2^{j-1}]$
  does the work for sufficiently large $n$. 
	
Recall that a Furstenberg family $\A$ is called \emph{partition regular} if given $A\in \A$ and any finite partition of $A$ as $A=\cup_{j=0}^m A_j$ implies there exists some $j$ with $0\leq j \leq m$ such that $A_j\in \A$. It is well known that $\overline{\mathcal{D}}, \overline{\mathcal{BD}}$ and $\PS$ are partition regular, see \cite{HiSt}. The family $\A$ is called \emph{shift invariant} if for every $i\in \Z$ and each $A\in \A$, we have $(A-i)\cap \Z_+\in \A$.
	
		Let $X$ be a separable metrizable and complete topological vector space and denote by $\mathscr{B}(X)$ the set of all bounded and linear operators on $X$. From now on, any operator considered here is in $\mathscr{B}(X)$. We are concerned with the study of the dynamics of continuous and linear operators acting on $X$. If $T$ is such an operator on $X$, $T$ is called \emph{hypercyclic} if there exists some $x\in X$ such that the recurrence set $N(x, U):=\{n\geq 1: T^n x \in U\}\neq \emptyset$, for any non-empty open set $U$ of $X$. Such $x$ is called a hypercyclic vector for $T$, and the set of such points is denoted by $HC(T)$. By Birkhoff's Transitivity Theorem, $T$ is hypercyclic on $X$ if and only if it is \emph{transitive}, i.e. the recurrence set $N(U, V)=\{n\geq 1: T^n(U)\cap V\neq \emptyset\}$ is non-empty, for any pair of non-empty open sets $U$ and $V$ of $X$. Hypercyclicity has been the subject of consistent investigation in the last decades. We refer the reader to \cite{GrPe}, \cite{BaMa} for a rich source on this subject. It is natural to wonder about those operators whose recurrence sets belong to a specific Furstenberg family of subsets of $\Z_+$.
	
	In \cite{BMPP2}, the notion of $\A$-transitivity was studied for operators in $\mathscr{B}(X)$. An operator $T\in \mathscr{B}(X)$ is called \emph{$\A$-transitive} if for any non-empty open sets $U, V$ of $X$, the set $N(U, V):=\{n\geq 0: T^n(U)\cap V\neq \emptyset\}\in \A$.
	Very recently the authors of \cite{BMPP1} have introduced and studied the notion of $\A$-hypercyclicity. 

Analogously, an operator $T$ is called \emph{$\A$-hypercyclic} if there exists $x\in X$ such that $N(x, U)\in \A$ for any non-empty open set $U$ of $X$. Such $x$ is called an $\A$-hypercyclic vector for $T$, and the set of such points is denoted by $\A HC(T)$. This notion was introduced and studied in \cite{BMPP1}, and it is indeed a generalization of the \emph{frequent hypercyclicity} notion (i.e. $\underline{\mathcal{D}}$-hypercyclicity). Frequently hypercyclic operators were introduced by Bayart and Grivaux in 2006, see \cite{BaGr1},\cite{BaGr2}. Since then it has been the most investigated class among the $\A$-hypercyclic operators, and so much have been studied in its connection with ergodic theory, among the more recent work on this being \cite{GrMa} \cite{GrMaMe}. More recently, other instances of $\A$-hypercyclicity have been considered in the literature like \emph{$\mathfrak{U}$-frequent hypercyclicity} (i.e. $\overline{\mathcal{D}}$-hypercyclicity) and \emph{reiterative hypercyclicity} (i.e. $\overline{\mathcal{BD}}$-hypercyclicity). Please refer to \cite{Sh}, \cite{BMPP1}, \cite{Me}, \cite{BoGr} for more information. Clearly, frequently hypercyclic operators are $\mathfrak{U}$-frequently hypercyclic, and these in turn are reiteratively hypercyclic. Any of these implications do not hold on the reverse direction, and the first counterexamples are rather intricate in their construction and with no apparent connection between them, see \cite{BaRu}, \cite{BMPP1}; however, very recently Bonilla and Grosse-Erdmann \cite{BoGr} exhibited a unified and much simpler way of obtaining such counterexamples.

	Recall that a point $x\in X$ is said to be \emph{periodic} for $T\in \mathscr{B}(X)$ if there exists $k\geq 1$ such that $T^kx=x$. An operator $T\in \mathscr{B}(X)$ is said to be \emph{chaotic} if it is hypercyclic and has a dense sets of periodic points.
	We recall to the reader that in \cite{Me} is proven that any chaotic operator is reiteratively hypercyclic. Moreover, a closer look at the proof reveals that in fact any chaotic operator is $\PS$-hypercyclic. And obviously, any $\PS$-hypercyclic operator is reiteratively hypercyclic. By Bayart and Grivaux \cite{BaGr2} there exists a frequently hypercyclic operator on $c_0(\Z_+)$ which is not chaotic. And more recently, one of the most important problems in linear dynamics has been solved by Menet \cite{Me} by showing the existence of a chaotic operator on $l^1(\Z_+)$ which is not $\mathfrak{U}$-frequently hypercyclic (hence no frequently hypercyclic). 
  
	Weighted shifts are a particular class of linear operators with a significant importance in linear dynamics. Each bilateral bounded weight $w=(w_k)_{k\in \Z}$ induces a \emph{bilateral weighted backward shift} $B_w$ on $X=c_0(\Z)$ or $\ell^p(\Z) (1\leq p<\infty)$, given by $B_{w}e_k:=w_{k}e_{k-1}$, where $(e_k)_{k\in \Z}$ denotes the canonical basis of $X$.
	Similarly, each unilateral bounded weight $w=(w_n)_{n\in \Z_+}$ induces a \emph{unilateral weighted backward shift} $B_w$ on $X=c_0(\Z_+)$ or $\ell^p(\Z_+) (1\leq p<\infty)$, given by $B_{w}e_n:=w_{n}e_{n-1}, n\geq 1$ with $B_{w}e_0:=0$, where $(e_n)_{n\in \Z_+}$ denotes the canonical basis of $X$.

  
  \section{Piecewise syndeticity level of recurrence sets of $\A$-hypercyclic operators}
	
	According to \cite{BMPP2}, if $T$ is an $\A$-transitive operator then $T$ is thick-transitive if and only if for any non-empty open sets $U$ and $V$ and any $N\geq 1$, there exists $B\in \A$ such that $(B+[-N, N])\cap \Z_+\subseteq N(U, V)$. So, for these operators the thickness level of any set $N(U, V)$ turns out to be high in such a way that the family $\A$ is also involved. On the other hand, as showed in Proposition 3 \cite{BMPP1}, no $F$-space $X$ supports a thick-hypercyclic operator. So, we cannot hope to have an operator $T\in \mathscr{B}(X)$ such that for some $x\in X$, the set $N(x, U)$ is thick for any non-empty open set $U$ of $X$. However, there exist piecewise syndetic hypercyclic operators, for instance any chaotic operator is piecewise syndetic hypercyclic, and the proof of Theorem 1.1 \cite{Me} essentially shows this fact. So, we could wonder how the piecewise syndeticity level of recurrence sets $N(x, U)$ of operators being both $\A$-hypercyclic and piecewise syndetic hypercyclic is, and eventually how the family $\A$ affects this level. In case this happens, we will be interested in how to take advantage of this fact.

Our aim in this section is to study the piecewise syndeticity level of recurrence sets of $\A$-hypercyclic operators, for a Furstenberg family $\A$ on $\Z_+$. In order to do so, we need to introduce the following family, that will  be an indicator of the piecewise syndeticity level of the family $\A$. So, define \[\PS^\A:=\{A\subseteq \N: \exists b \quad \forall L \quad \{z \geq 1: [z, z+L)\subset \cup_{t=1}^b A-t\}\in \A\}.\]

In general, a set $A\in \PS\cap \A$ not necessarily is in $\PS^{\A}$. However, for recurrence sets of $\A$-hypercyclic operators this is not the case as shown below in Proposition \ref{high_PS_level}. 

Furthermore, suppose $T$ is both frequently hypercyclic and $\PS$-hypercyclic. As mentioned in Corollary \ref{charact_PS_Hyp} below, the set of $\PS$-hypercyclic vectors contains a dense $G_\delta$-set. On the other hand, the set of frequently hypercyclic vectors is known to be meager in $X$ \cite{BaRu}. Thus, it might happen that no $\PS$-hypercyclic vector is frequently hypercyclic. This situation must not happen as the next result shows.

\begin{proposition}
\label{high_PS_level}
Let $\A$ be a Furstenberg family on $\Z_+$ and $T\in \mathscr{B}(X)$. 

(a) Consider the following:

(i) $T$ is both $\A$-hypercyclic and $\PS$-hypercyclic,

(ii) $T$ is $\A$-hypercyclic and for any non-empty open set $U$ of $X$ there exists $y\in X$ such that $N(y, U)\in \PS$,

(iii) $T$ is $\PS^{\A}$-hypercyclic.

Then, (i) $\Rightarrow$ (ii) $\Rightarrow$ (iii), and the set of $\A$ hypercyclic vectors for $T$ is contained in the set of $\PS^\A$ hypercyclic vectors for $T$.

(b) Moreover, if $\A$ is partition regular and shift invariant, then (i)-(iii) are equivalent, and the set of $\A$ hypercyclic vectors for $T$ coincides with the set of $\PS^\A$ hypercyclic vectors for $T$.

\end{proposition}

\begin{proof}
Obviously $(i)$ implies $(ii)$. We just need to show that $(ii)$ implies $(iii)$. Let $x$ be a $\A$ hypercyclic vector for $T$ and $U$ a non-empty open set of $X$. Then there exists $y\in X$ such that $N(y, U)\in b$-$\PS$ for some $b\geq 1$, i.e. for any $L\geq 1$ there exists $z\geq 1$ such that $[z, z+L)\subset N(y, \cup_{t=1}^b T^{-t}U)$. Set $U_L=\cap_{l=0}^{L-1} \cup_{t=1}^b T^{-(l+t)}U$, it is a non-empty open set of $X$, since $T$ is continuous and $T^zy\in U_L$. Since $x$ is $\A$ hypercyclic vector for $T$, we get $A_L:=N(x, U_L)\in \A$. In other words, $N(x, U)\in \PS^{\A}$. 
\end{proof}

If we consider $\A$ the family of infinite subsets of $\N$ in Proposition \ref{high_PS_level}, we then automatically have.

\begin{corollary}
\label{charact_PS_Hyp}
An operator $T\in \mathscr{B}(X)$ is $\PS$-hypercyclic if and only if $T$ is hypercyclic and  for any non-empty open set $U$ of $X$ there exists $y\in X$ such that $N(y, U)\in \PS$. Furthermore, the set of $\PS$-hypercyclic vectors is a dense $G_\delta$-set of $X$, as it coincides with the set of hypercyclic vectors for $T$.
\end{corollary}

\begin{corollary}
Let $T\in \mathscr{B}(X)$ be a $\PS$-hypercyclic operator, $\A$ a Furstenberg family on $\Z_+$ and $k\geq 1$. Consider the following:

(i) the $k$-fold product operator $T\times \cdots \times T$ is $\A$-hypercyclic,

(ii) the $k$-fold product operator $T\times \cdots \times T$ is $\PS^{\A}$-hypercyclic.

Then (i) implies (ii). Moreover, if $\A$ is partition regular and shift invariant, then (i)-(ii) are equivalent.

\end{corollary}
\begin{proof}
We just need to show $(i)$ implies $(ii)$. The case $k=1$ is just Proposition \ref{high_PS_level}. Suppose $k=2$, then for any non-empty open sets $U, V$ there exist $x, y\in X$ such that $\{n\geq 0: T^n \times T^n (x, y)\in U\times V\}\in \PS$. Indeed, pick $m\in N(V, U)$ and define the non-empty open set $V_m=T^{-m}V\cap U$. Hence for any $n\in N(x, V_m)\in \PS$ we have $T^n x\in U$ and $T^n y\in V$ for $T^m x=y$. On the other hand, $T\times T$ is $\A$-hypercyclic thus by Proposition \ref{high_PS_level} we get $T\times T$ is $\PS^{\A}$-hypercyclic. If $k>2$, repeating the same argument we obtain that the $k$-fold product $T\times \cdots \times T$ is $\PS^{\A}$-hypercyclic.
\end{proof}

 \subsection{Proof of Theorem \ref{theorem.irregular.BD}}
As already mentioned, the proof of Theorem \ref{theorem.irregular.BD} is combinatorial number theory in essence. Indeed, the main ingredient of the proof is a convenient variation of a result given by Hindman in 1990 which follows our purposes. In \cite{Hi}, Hindman showed that if $A\subseteq \N$ is such that $\alpha=\underline{d}(A)\leq \overline{Bd}(A)=\gamma$ then for any $\epsilon >0$ there exists $b\geq 1$ such that $\underline{d}(\cup_{t=1}^b A-t)\geq \alpha/\gamma -\epsilon$. It is essential to us to be able to get rid of the dependence of $b$ on $\epsilon$, in order to get a result in the vein of Theorem \ref{theorem.irregular.BD}. We show in the next lemma that by adding the hypothesis $A\in \PS$ we got $b$ not depending on $\epsilon$ in Hindman's result, which is exactly what we need. At this point it might then be clear our interest in introducing and studying the notion of $\PS$-hypercyclicity.

\begin{lemma}
\label{lemma_cefalu}
 Let $A\in \PS$  with $0<\alpha=\underline{d}(A)\leq \overline{Bd}(A)=\gamma$. Then there exists $b' \geq 1$ such that $\overline{d}(\cup_{t=1}^{b'} A-t)\geq \alpha/\gamma$.
  \end{lemma}
\begin{proof}

 Let $A\in \PS$, then there exists $b\geq 1$ such that $\cup_{t=1}^b A-t$ is thick. Set $B=\{z\in A:[z-b, z+b]\cap A=\emptyset\}$ and $\hat{A}=A\setminus B$. 

Hence, there exist unique sequences $(z_i)_i$ and $(r_i)_i$ with $(z_i)$ strictly increasing and $r_i\geq 1$ such that

\begin{itemize}
\item $ \cup_{t=1}^b \hat{A}-t=\cup_{i\geq 1}[z_i, z_i+r_i),$

\item $z_{i+1}-(z_i+r_i)>b  \quad \forall i\geq 1,$

\item for any $N\geq 1$, there exists $j\geq 1$ such that $r_j\geq N.$
\end{itemize}

By hypothesis $\overline{Bd}(A)=\gamma$. Let $\zeta>0$, we have $\overline{Bd}(\hat{A})\leq \overline{Bd}(A)<\gamma+\zeta$, which means that by (\ref{condition_etna}) we can pick $m'\geq 1$ such that for any $z\geq 1$ and any $r\geq m'$ we have 
\begin{equation}
\label{condition_modica}
\abs{\hat{A}\cap [z, z+r)}<(\gamma+\zeta) \cdot r.
\end{equation}

On the other hand, by hypothesis $\underline{d}(A)=\alpha> \alpha-\zeta$.
Hence, there exists $m\geq 1$ such that for any $n\geq m$, 
\begin{equation}
\label{condition_scicli}
\abs{A\cap [1, n]}> (\alpha-\zeta)\cdot n.
\end{equation}

Set $I:=\{i\geq 1: r_i\geq m'\}$ and define $\hat{A}_1\subseteq \hat{A}$ such that $ \cup_{t=1}^b \hat{A}_1-t=\cup_{i\in I}[z_i, z_i+r_i)$.

According to the properties of $(r_i)$ listed above, for any $M\geq m$ there exist $z$ and $\nu$ with $z\geq M$ and $\nu\geq z$ satisfying the following both conditions:

\begin{itemize}

\item $\abs{\big(\cup_{t=1}^{2m'b} A-t\big) \cap [1, z]}\geq \abs{\big(\cup_{t=1}^{b} (\hat{A}_1\cup B)-t\big) \cap [1, \nu]},$

\item $\abs{A\cap [1, z]}\leq \abs{(\hat{A}_1\cup B)\cap [1, \nu]}$.
\end{itemize}

Set $j=\max\{i:z_i+r_i< \nu\}$. Note that $\nu$ can be chosen such that $\nu-z_{j+1}\geq m'b$.

Thus, combining (\ref{condition_modica}), (\ref{condition_scicli}) and considering that $1/b\leq \gamma+\zeta$ we have 
\[\abs{(\cup_{t=1}^{2m'b} A-t)\cap [1, z]}\geq \abs{\big(\cup_{t=1}^{b} (\hat{A}_1\cup B)-t\big) \cap [1, \nu]} \geq \]

\[\abs{\cup_{t=1}^b\big((\hat{A}_1\cup B)-t\big)\cap ([1, \nu]\setminus \cup_{i=1}^j [z_i, z_i+r_i))} + \sum_{i\in I \cap [1, j]} r_i + (\nu-z_{j+1})= \]

\[ b\abs{B\cap [1, \nu]}  +  \sum_{i\in I \cap [1, j]} r_i + (\nu-z_{j+1})>\]

\[b\abs{B\cap [1, \nu]} + \frac{\sum_{i\in I \cap [1, j]} \abs{\hat{A}_1\cap [z_i, z_i+r_i)}+\abs{\hat{A}_1\cap [z_{j+1}, \nu]}}{\gamma+\zeta} \geq\]

\[\frac{\abs{B\cap [1, \nu]} + \sum_{i\in I \cap [1, j]} \abs{ \hat{A}_1\cap[z_i, z_i+r_i)}+\abs{\hat{A}_1\cap [z_{j+1}, \nu]}}{\gamma+\zeta}=\]

\[\frac{\abs{(\hat{A}_1\cup B)\cap[1, \nu]}}{\gamma+\zeta}\geq \frac{\abs{ A\cap[1, z]}}{\gamma+\zeta} > \frac{(\alpha-\zeta)}{(\gamma+\zeta)}\cdot z.\]

 Set $b':=2m'b$. This means $\overline{d}(\cup_{t=1}^{b'} A-t)> \frac{\alpha-\zeta}{\gamma+\zeta}$. Letting $\zeta$ tend to $0$ we get $\overline{d}(\cup_{t=1}^{b'} A-t)\geq \alpha/\gamma$.
\end{proof}

\begin{remark}
		In Theorem 3.2 \cite{Hi}, Hindman showed that given $0<\alpha\leq\gamma\leq 1$, there exists $A\subseteq \N$ with $\underline{d}(A)=\overline{d}(A)=\alpha, \overline{Bd}(A)=\gamma$ and for all $b\in \N, \overline{d}(\cup_{t=1}^b A-t)\leq \alpha / \gamma$, with a strict inequality if $\gamma <1$. So, Lemma \ref{lemma_cefalu} does not necessarily hold if we drop the condition about the piecewise syndeticity of $A$.
		\end{remark}

	Now we have all we need to prove Theorem \ref{theorem.irregular.BD}.

\textit{Proof of Theorem \ref{theorem.irregular.BD}.}
Let $T, x$ and $U$ as in the hypothesis of Theorem \ref{theorem.irregular.BD}. Suppose $x$ is a frequently hypercyclic vector for $T$. Since $T$ is $\PS$-hypercyclic, by Proposition \ref{high_PS_level} we have that $N(x, U)\in \PS$. Thus, $N(x, U)\in \underline{\mathcal{D}}\cap \PS$. Suppose that $\underline{d}(N(x, U))=\overline{Bd}(N(x, U))>0$, then by Lemma \ref{lemma_cefalu}  there exists $b\geq 1$ such that $\overline{d}(\{k\geq 1: T^k x\in \cup_{t=1}^b T^{-t}U\})=1$. Set $W=\cup_{t=1}^b T^{-t}U$, then $W$ is non-dense on $X$ by hyphotesis. So, we can pick a non-empty open set $V$ such that $V\cap W=\emptyset$. Then $\underline{d}(N(x, V))=0$ which contradicts the fact that $x$ is a frequently hypercyclic vector. Hence, $0<\underline{d}(N(x, U))<\overline{Bd}(N(x, U))$. Now, suppose $x$ is a hypercyclic vector but not frequently hypercyclic vector, then by Corollary \ref{charact_PS_Hyp}, $x$ is $\PS$-hypercyclic vector, so $0=\underline{d}(N(x, U))<\overline{Bd}(N(x, U))>0$. This completes the proof of  Theorem \ref{theorem.irregular.BD}.

 We would like to remark that there exists a frequently hypercyclic and $\PS$ hypercyclic weighted backward shift on $c_0(\Z)$ satisfying the conclusion of Theorem \ref{theorem.irregular.BD} for any non-empty open set $U$ of $c_0(\Z)$.

First, we need the following. In Theorem 12 \cite{BaRu} Bayart and Ruzsa give a characterization for frequently hypercyclic and $\mathfrak{U}$-frequently hypercyclic weighted backward shifts on $c_0(\Z)$. Taking into consideration that $\PS$ is partition regular, following the same lines of the proof of Theorem 12 \cite{BaRu} we obtain a characterization for $\PS$-hypercyclic weighted backward shifts on $c_0(\Z)$, which reads as follows.

\begin{proposition}
\label{char.PS.hyp}
Let $w=(w_n)_{n\in \Z}$ be a bounded and bounded below sequence of positive integers. Then $B_w$ is $\PS$-hypercyclic if and only if there exist a sequence $(M(p))$ of positive real numbers tending to $+\infty$ and a sequence $(E_p)$ of subsets of $\Z_+$ such that 

(a) For any $p\geq 1, E_p\in \PS$;

(b) For any $p, q\geq 1, p\neq q, (E_p+[-p, p]) \cap (E_q+[-q, q])=\emptyset$;

(c) $\lim_{n\to +\infty, n\in E_p} w_1\cdots w_n=+\infty$;

(d) For any $p, q \geq 1$, for any $n\in E_p$ and any $m\in E_q$ with $n\neq m$,

\[ \begin{cases} 
      w_1\cdots w_{m-n}\geq M(p)M(q) & \mbox{ provided } m>n \\
      w_{m-n+1}\cdots w_0\leq \frac{1}{M(p)M(q)} & \mbox{ provided } m<n. \\
    
   \end{cases}
\]
\end{proposition}

Then we have the following.

\begin{proposition}
\label{riverside}
There exists a frequently hypercyclic and $\PS$ hypercyclic weighted backward shift $B_w$ on $c_0(\Z)$, for which $\underline{d}(N(x, U))<\overline{Bd}(N(x, U))$ for any hypercyclic vector $x$ and any non-empty open set $U$ of $c_0(\Z)$.
\end{proposition}

\begin{proof}

	 Indeed, the operator constructed by Bayart and Ruzsa in Theorem 7 \cite{BaRu} is the one we are considering here. Bayart and Ruzsa exhibited there a weight $w$ and  a sequence of sets $(E_p)$ defined as follows: 
\[
E_p=\bigcup_{u\in A_p} I^{a, \epsilon}_u \cap b_p \N;
\]
where each $A_p$ is syndetic with $\bigcup_{p\geq 1}A_p$ a partition of $\N$, $(b_p)$ is some increasing sequence of positive numbers and for some $a>1$ and $\epsilon >0, I^{a, \epsilon}_u:=[(1-\epsilon)a^u, (1+\epsilon)a^u]$.

They showed that the corresponding weighted backward shift $B_w$ is frequently hypercyclic and the sets $E_p$ satisfy conditions (b), (c) and (d).           

 Note that $E_p\in \PS$, for any $p\geq 1$. Hence, by Proposition \ref{riverside} we conclude that $B_w$ is also $\PS$-hypercyclic.

 Furthermore, Grivaux and Matheron pointed out in Section 4.1 \cite{GrMa} that this operator also satisfies $\sup_{R>0}\overline{d}(N(x, B_R))<1$, for any hypercyclic vector $x$; where $B_R$ denotes the closed ball with radius $R$ and centered at $0$. Then, we have that
\begin{equation}
\label{cond_BaRu}
\sup_{s\geq 1}\overline{d}(N(x, \cup_{t=1}^sT^{-t}U))<1,
\end{equation}
 for any hypercyclic vector $x$ and any non-empty open set $U$ of $c_0(\Z)$. Now, let $x$ be a frequently hypercyclic vector of $B_w$ and $U$ a non-empty open set of $c_0(\Z)$. Suppose that $\underline{d}(N(x, U))=\overline{Bd}(N(x, U))>0$, then by  Proposition \ref{high_PS_level} and Lemma \ref{lemma_cefalu}  there exists $b\geq 1$ such that $\underline{d}(\cup_{t=1}^b N(x, U)-t)=\underline{d}(\{k\geq 1: T^k x\in \cup_{t=1}^b T^{-t}U\})=1$, which contradicts (\ref{cond_BaRu}). Hence, $0<\underline{d}(N(x, U))<\overline{Bd}(N(x, U))$. Now, suppose $x$ is a hypercyclic vector but not frequently hypercyclic vector, then by Corollary \ref{charact_PS_Hyp}, $x$ is $\PS$-hypercyclic vector, so $0=\underline{d}(N(x, U))<\overline{Bd}(N(x, U))>0$.

\end{proof}

	
  \subsection{How to be piecewise syndetic hypercyclic}
	Due to Theorem \ref{theorem.irregular.BD} it is natural to be interested to know when an operator is piecewise syndetic hypercyclic. In this respect we have Corollary \ref{charact_PS_Hyp}. Our aim in this subsection is to try to understand better the conclusion of Corollary \ref{charact_PS_Hyp} by getting variations of this result.

Thanks to Baire Category Theorem, Bonilla and Grosse-Erdmann \cite{BoGr} have recently showed that $T$ is hypercyclic if and only if for any non-empty open set $V$ and any non-empty open set $O$, there exists $x\in O$ such that $N(x, V)$ is non-empty.

 This leads us to the following question: in order to get $\PS$-hypercyclicity, how  Corollary \ref{charact_PS_Hyp} get affected if we replace the condition \emph{$T$ hypercyclic} by the weaker condition \emph{for any open neighborhood $W$ of 0 and any non-empty open set $O$ of $X$, there exists $x\in O$ such that $N(x, W)$ is non-empty}?

The following result gives us the answer.
	
	\begin{proposition} 
	\label{prop.padova}
	(a) Let $T\in \mathscr{B}(X)$ such that for any non-empty open set $V$ of $X$ there exists $y\in X$ such that $N(y, V)$ is syndetic. If for any open neighbourhood $W$ of $0$ and any non-empty open subset $O$ of $X$, there exists $x\in O$ such that $N(x, W)$ is non-empty, then $T$ is $\PS$-hypercyclic;
	
	(b) Let $T\in \mathscr{B}(X)$ such that for any open neighbourhood $W$ of $0$ and any non-empty open set $O$ of $X$ there exists $x\in O$ such that $N(x, W)$ is syndetic. If for any non-empty open set $V$ of $X$ there exists $y\in X$ such that $N(y, V)\in \PS$, then $T$ is $\PS$-hypercyclic.
	\end{proposition}
	
	\begin{proof}
	We are following the same idea of the proof of Theorem 5.16 \cite{GrMaMe}. We give all details for the sake of completeness. First, let us denote by $\mathfrak{I}_N$ the set of intervals on $\N$ of length bigger than $N$. 
	
	  By a Baire Category argument, it is enough to show that for any non-empty open set $V$ of $X$ there exists a natural number $b_V\geq 1$ such that for any $N\geq 1$ the set 
	\[ G_{V, N}:=\{u\in X: \exists J\in \mathfrak{I}_N,  J\subseteq \cup_{t=0}^{b_V}(N(u, V)-t)\}
	\]
	is dense in $X$. In fact, if $(V_q)_{q\geq 1}$ is a countable basis of open sets of $X$, the set $\cap_{q, N}G_{V_q, N}$ is a dense $G_{\delta}$ set consisting of $\PS$-hypercyclic vectors for $T$. So, let $V$ be a non-empty open set, $z\in V$ and $\epsilon >0$ such that $B(z, \epsilon)\subseteq V$, where $B(z, \epsilon)$ denotes the open ball centered at $z$ with radius $\epsilon$. 
	
	Let us show first part (a). By hypothesis, there exists $y\in X$ such that $N(y, B(z, \epsilon/2))$ is syndetic with gap bounded by $b_V$.
	
	Fix $N\geq 1$ and $U$ a non-empty open set of $X$, we need to show that $G_{V, N}\cap U\neq\emptyset$.
	Set $O=U-y$. By assumption, there exists $x\in O$ such that $N(x, B(0, \epsilon/2))$ is non-empty, hence thick. Indeed, this follows from the fact that $B(0, \epsilon/2)$ is an open neighbourhood of $0$ and the continuity of the operator $T$.
	Then $u=x+y\in U$. Moreover, $u\in G_{V, N}$. Indeed, the set $N(y, B(z, \epsilon/2)) \cap N(x, B(0, \epsilon/2))\in \PS$ with gap bounded by $b_V$ since any thick set intersects any syndetic set with gap bounded by $b$ in a $b$-$\PS$ set. Hence, there exists an interval $J\in \mathfrak{I}_N$ such that 
	\[
	J\subseteq \cup_{t=0}^{b_V}\Big(\big\{n\in J: \norm{T^nu-z}\leq \norm{T^nx}+\norm{T^ny-z}\leq \epsilon\big\}-t\Big).
	\]
	
 To establish part (b), by hypothesis there exists $y\in X$ such that the set 
 $N(y, B(z,\epsilon/2))$ is piecewise syndetic with gap bounded by $b_V$. Fix $N\geq 1$ and $U$ a non-empty open set of $X$, we need to show that $G_{V, N}\cap U\neq\emptyset$.
	Set $O=U-y$. By assumption, there exists $x\in O$ such that $N(x, B(0, \epsilon/2))$ is syndetic, hence for every $N$ there exists a syndetic set $B_N$ such that $B_N+[-N, N]\subseteq N(x, B(0, \epsilon/2))$. Indeed, this follows from the fact that $B(0, \epsilon/2)$ is an open neighbourhood of $0$ and the continuity of the operator $T$.
	 Note that the set $N(y, B(z, \epsilon/2)) \cap N(x, B(0, \epsilon/2))\in \PS$ with gap bounded by $b_V$. In fact, let $M\geq N$ then there exists a syndetic set $B_M$ with gap bounded by $b_M$ such that $B_M+[0, M]\subseteq N(x, B(0, \epsilon/2))$. We can pick $r\geq 1$ such that $[r, r+b_M+M]\subseteq \cup_{t=0}^{b_V} N(y, B(z, \epsilon/2))-t$ since $N(y, B(z, \epsilon/2))$ is $b_V$-$\PS$. Hence, there exists some $j\in [0, b_M]$ such that 
	\[
	[r+j, r+j+M]\subseteq \big(\cup_{t=0}^{b_V} N(y, B(z, \epsilon/2))-t\big)\cap N(x, B(0, \epsilon/2)).
	\]
	
	Letting $u=x+y$, obviously $u\in U$. Moreover, $u\in G_{V, N}$ since there exists an interval $J\in \mathfrak{I}_N$ such that 
	\[
	J\subseteq \cup_{t=0}^{b_V}\Big(\big\{n\in J: \norm{T^nu-z}\leq \norm{T^nx}+\norm{T^ny-z}\leq \epsilon\big\}-t\Big).
	\]

	\end{proof}
	As a consequence we can strengthen the conclusion of Thereom 5.16 \cite{GrMaMe}.
	\begin{corollary} 
	\label{prop.verona}
	Let $T\in \mathscr{B}(X)$ such that for any non-empty open set $V$ of $X$ there exists $y\in X$ such that $N(y, V)$ is syndetic. Assume there exists $\alpha>0$ such that for any open neighbourhood $W$ of $0$ and any non-empty open subset $O$ of $X$, one can find $x\in O$ and an arbitrarily large integer $n$ such that $[n, n+\alpha n]\subset N(x, W)$. Then $T$ is $\mathfrak{U}$-frequently hypercyclic and $\PS$-hypercyclic. In particular, $T$ is  $\PS^{\overline{\mathcal{D}}}$-hypercyclic.
	\end{corollary}
	\begin{proof}
	By Proposition \ref{high_PS_level} it is enough to show that $T$ is both $\mathfrak{U}$-frequently hypercyclic and $\PS$-hypercyclic. Note that in particular, $T$ satisfies hypothesis of Proposition \ref{prop.padova}, so $T$ is $\PS$-hypercyclic. To establish that $T$ is $\mathfrak{U}$-frequently hypercyclic is essentially the content of Theorem 5.16 \cite{GrMaMe}, and we refer the reader to its original proof. Just note that we have just weakened a little bit the first hypothesis in Theorem 5.16 \cite{GrMaMe}.
	\end{proof}
		
	Finally, we would like to point out to the reader a sufficient condition for an operator to be frequently hypercyclic and chaotic (hence frequently hypercyclic and $\PS$-hypercyclic) obtained in Theorem 5.31 \cite{GrMaMe}.
	
		\section*{Aknowledgements}
		
		We would like to thank Michael Lin, Tom Meyerovitch and Guy Cohen for their hospitality received during the preparation of this paper and to Sophie Grivaux for pointing us a mistake in an earlier version of Lemma \ref{lemma_cefalu}.

\end{document}